\documentclass{amsart}

\newcommand\myref[1]{\cite{#1}}

\newtheorem{theorem}{Theorem}[section]
\newtheorem{proposition}[theorem]{Proposition}

\newtheorem{lemma}[theorem]{Lemma}

\theoremstyle{remark}

\newtheorem{remark}[theorem]{Remark}

\newcommand\NS{\mathop{\rm NS} \nolimits}
\newcommand\Pic{\mathop{\rm Pic} \nolimits}
\newcommand\Spec{\mathop{\rm Spec} \nolimits}

\newcommand\rk{\mathop{\rm rk} \nolimits}
\newcommand\Tr{\mathop{\rm Tr} \nolimits}
\newcommand\disc{\mathop{\rm disc} \nolimits}

\newcommand\Q{\mathbb{Q}}
\newcommand\Qbar{\overline{\mathbb{Q}}}
\newcommand\D{\mathcal{D}}
\newcommand\R{\mathbb{R}}

\newcommand\F{\mathbb{F}}
\newcommand\A{\mathbb{A}}
\newcommand\p{\mathfrak{p}}
\renewcommand\P{\mathbb{P}}
\newcommand\Xbar{\overline{X}}
\renewcommand\O{\mathcal{O}}
\newcommand\kbar{\overline{k}}
\newcommand\Fbar{\overline{\F}}
\newcommand\Ybar{\overline{Y}}
\newcommand\Xint{\mathfrak{X}}

\newcommand\et{\text{\rm \'et}}
\newcommand\ra{\rightarrow}
\renewcommand\vec{\mathbf}
\newcommand\error{0.1}

\def\ww{\omega}
\def\ss{\sigma}
\def\cc{\gamma}
\def\Bbb#1{{\mathbb #1}}

\def\mymatrix#1{\begin{matrix}#1\end{matrix}}

\begin{document}

\title[Canonical vector heights on K3 surfaces]
      {Canonical vector heights on K3 surfaces with Picard number
       three -- addendum}
 \author{Arthur Baragar}
 \author{Ronald van Luijk}
 \address{University of Nevada Las Vegas, Las Vegas, NV 89154-4020}
  \email{baragar@unlv.nevada.edu}
 \address{Mathematical Sciences Research Institute,
          17 Gauss Way, Berkeley, CA 94720-5070}
  \email{rmluijk@msri.org}
 \thanks{The first author is supported in part by NSF grant
 DMS-0403686.}
 \copyrightinfo{2006}{American Mathematical Society}

\begin{abstract}
In an earlier paper by the first author, an argument for the
nonexistence of canonical vector heights on K3 surfaces of Picard
number three was given, based on an explicit surface that was not
proved to have Picard number three. In this paper, we fill the gap
in the argument by redoing the computations for another explicit
surface for which we prove that the Picard number equals three.  The
conclusion
remains unchanged.  
\end{abstract}

\maketitle

\section{Introduction}

In \cite{barorig} the first author gave convincing numerical
evidence for the nonexistence of canonical vector heights on K3
surfaces of Picard number $3$. The intent of this paper is to fill a
gap in the argument, which was pointed out by Yuri Tschinkel in the
review of the paper, and privately by Bert van Geemen.

As in \cite{barorig}, the Picard number of a surface will always
mean the geometric Picard number. The Picard number of the explicit
K3 surface $V$ used in \cite{barorig} is at least $3$, but was not
proved to equal $3$. Since $3$ is odd, the only currently known
method to prove that this lower bound is sharp requires two primes
of good reduction for $V$, such that the Picard number of both
reductions equals $4$, see \cite{picone}. Modulo $2$ and $3$ the
Picard numbers turn out to be $16$ and $6$ respectively (depending
on Tate's conjecture).  The computations required to calculate the
Picard number modulo larger primes are currently beyond our ability.
    This is
why in the next section we construct a new example $Y$ for which we
can use the primes $2$ and $3$ to prove that the Picard number
equals $3$. In the last section we redo the necessary computations
for this example $Y$, referring to \cite{barorig} for details. We
compute various canonical heights, which we believe to be correct up
to an error of at most $0.0001$. Our main theorem states that if the
errors are at most $\error$, then a canonical vector height on $Y$
does not exist. This also suggests that, except perhaps in very
special cases%
, a K3 surface with Picard number at least three will not admit a
canonical vector height.

\section{A K3 surface with Picard number three}

Let $k$ be a field with a fixed algebraic closure $\kbar$. Let
$X$ be a smooth surface over $k$ in $\P^1 \times \P^1 \times \P^1$, given by
a $(2,2,2)$-form. Then $X$ is a K3 surface, which implies that linear,
algebraic, and numerical equivalence all coincide. This means that the Picard
group $\Pic \Xbar$ and the N\'eron-Severi group $\NS \Xbar$ of $\Xbar
= X_{\kbar}$ are naturally isomorphic, finitely generated, and
free. Their rank is called the (geometric) Picard number of $X$.
By the Hodge Index Theorem, the intersection pairing gives this
group the structure of a lattice with signature $(1, \rk \NS \Xbar -1)$.
For detailed definitions of all these notions, see \cite{picone}.

For $i=1,2,3$, let $\pi_i \colon \, X \ra \P^1$ be the projection from $X$
to the $i$-th copy of $\P^1$ in $\P^1 \times \P^1 \times \P^1$.
%
%
Let $D_i$ denote the divisor class represented by a fiber of
$\pi_i$. We have $D_i \cdot D_j = 2$ for $i \neq j$ and
since any two different fibers of $\pi_i$ are disjoint,
we find
$D_i^2=0$. It follows that the intersection matrix $(D_i \cdot
D_j)_{i,j}$ has rank $3$, so the $D_i$ generate a subgroup of the
N\'eron-Severi group $\NS(\Xbar)$ of rank $3$. Our goal is to
find an explicit example where the rank of $\NS(\Xbar)$ equals $3$.

Let $x$, $y$, and $z$ denote the affine coordinates of $\A^1$ inside
the three copies of $\P^1$ in $\P^1 \times \P^1 \times \P^1$.
Let $Y/\Q$ be the surface given by $G_1x^2+G_2x+3G_3-2L_1L_2=0$ with
\begin{align*}
 G_1 & = -y^2z^2 + 3y^2z + 2y^2 - 2yz^2 + 3yz + 3y + 2z^2 + 2z - 1,
 \\
 G_2 & = 2y^2z^2 + 3y^2z + 3y^2 + 2yz^2 + 2yz + 3z^2 + z + 2, \\
 G_3 & = y^2z + y^2 + y + z^2 + z, \\
 L_1 & = yz - y - z, \\
 L_2 & = yz+1. \\
\end{align*}
\begin{theorem}\label{rankthree}
The surface $Y$ is smooth.
The Picard number of $Y_{\Qbar}$ equals $3$.
\end{theorem}

To bound the Picard number of $Y$ we use the method described in
\cite{picone}. We first state some results and notation that we will
use. Let $X$ be any smooth surface over a number field $K$ and let
$\p$ be a prime of good reduction with residue field $k$. Let
$\Xint$ be an integral model for $X$ over the localization $\O_\p$
of the ring of integers $\O$ of $K$ at $\p$. Let $k'$ be any
extension field of $k$. Then by abuse of notation we will write
$X_{k'}$ for $\Xint \times_{\Spec \O_\p} \Spec k'$.

\begin{proposition}\label{boundNS}
Let $X$ be a smooth surface over a number field $K$ and let $\p$
be a prime of good reduction with residue field $k$.
Let $l$ be a prime not dividing $q = \#k$. Let $F$ denote the
automorphism on $H^2_{\et}(X_{\kbar},\Q_l)(1)$ induced by $q$-th power
Frobenius. Then there are natural injections
$$
\NS(X_{\overline{K}}) \otimes \Q_l \hookrightarrow \NS(X_{\kbar})
\otimes \Q_l \hookrightarrow H^2_{\et}(X_{\kbar},\Q_l)(1),
$$
that respect the intersection pairing and the action of Frobenius
respectively. The rank of $\NS(X_{\kbar})$
is at most the number of eigenvalues of $F$ that are
roots of unity, counted with multiplicity.
\end{proposition}
\begin{proof}
See \myref{heron}, Prop. 6.2 and Cor. 6.4. Note that in the referred
corollary, Frobenius acts on the
cohomology group $H^2_{\et}(X_{\kbar},\Q_l)$ without a twist. Therefore,
the eigenvalues are scaled by a factor $q$.
\end{proof}

\begin{remark}\label{tate}
Tate's conjecture (see \cite{tate}) states that the rank of
$\NS(X_{\kbar})$ in Proposition \ref{boundNS}
is in fact equal to the number of eigenvalues
of $F$ that are roots of unity, counted with multiplicity.
\end{remark}

\begin{lemma}\label{uptosq}
If $\Lambda'$ is a sublattice of finite index in a lattice $\Lambda$,
then we have $\disc \Lambda' = [\Lambda : \Lambda']^2 \disc \Lambda$.
\end{lemma}

\proof[Proof of Theorem \ref{rankthree}]
 We write $Y_p$ and $\Ybar_p$ for $Y_{\F_p}$ and
$Y_{\Fbar_p}$ respectively.
One easily checks that $Y_p$ is smooth for $p=2$ and $p=3$, so
$Y$ itself is smooth and $Y$ has good reduction at $2$ and $3$.
%
%
Both $Y_2$ and $Y_3$ contain a fourth divisor class that is linearly
independent of the earlier described classes $D_i$ for $i=1,2,3$. On
$Y_2$ we have the curve $C_2$ parameterized by
$([x:1],[1:0],[1:1])$. On $Y_3$ we have the curve $C_3$ given by $x
= L_1 = 0$. For $p=2,3$, let $\Lambda_p$ denote the sublattice of
the N\'eron-Severi group of $\Ybar_p$ generated by $D_1,D_2,D_3$,
and $C_p$. The intersection matrices associated to the sequences of
classes $\{D_1,D_2,D_3,C_2\}$ and $\{D_1,D_2,D_3,C_3\}$ are
\[
\left[\begin{matrix} 0&2&2&1 \\
                     2&0&2&0 \\
                     2&2&0&0 \\
                     1&0&0&-2
    \end{matrix}\right]
\qquad \hbox{and} \qquad
    \left[ \begin{matrix} 0&2&2&0 \\
                          2&0&2&1 \\
                          2&2&0&1 \\
                          0&1&1&-2
    \end{matrix}\right],
\]
so $\Lambda_2$ and $\Lambda_3$ have discriminants $-28$ and $-32$
respectively. We will now show
that the Picard numbers of $\Ybar_2$ and $\Ybar_3$ both equal $4$.
Almost all fibers of the fibration $\pi_1$ are smooth curves of
genus $1$. Using {\sc magma} we counted the number of points over
small fields fiber by fiber. The total numbers of points are
given in Table \ref{tableone}.
\begin{table}[t]
$$
\begin{array}{|l|l|l|}
\hline
n  & \# Y_2(\F_{2^n}) & \# Y_3(\F_{3^n}) \cr
\hline
1  & 13   &17\cr
2  & 25   &107\cr
3  & 85   &848\cr
4  & 289  &6719\cr
5  & 1153 &60632\cr
6  & 4273 &536564\cr
7  & 16897&4793855\cr
8  & 65025&43091783\cr
9  & 266305&387501194\cr
10 & 1050625& \cr
\hline
\end{array}
$$
\caption{Number of points over some finite fields.} \label{tableone}
\end{table}
The Lefschetz Trace Formula relates the number of $\F_{p^n}$-rational
points on $Y_p$ to the traces of the $p^n$-th power Frobenius
acting on $H^i_\et(\Ybar_p,\Q_l)(1)$ for $i=0,\ldots,4$ by
$$
\#Y_p(\F_{p^n}) = \sum_{i=0}^4 (-p^{n/2})^i \cdot \left(
\mbox{trace of $p^n$-th power Frobenius on
  $H^i_\et(\Ybar_p,\Q_l)(1)$}\right).
$$
Normally this is phrased in terms of the cohomology groups without the twist.
For K3 surfaces we have $\dim H^i = 1,0,22,0,1$ for $i=0,1,2,3,4$
respectively.
Since the action for $i \neq 2$ is trivial,
from the numbers in Table \ref{tableone} we can compute the traces of
powers of the automorphism $F_p$ on $H^2_{\et}(\Ybar_p,\Q_l)(1)$ that is
induced by $p$-th power Frobenius.
We find $p^n \cdot \Tr F_p^n = \#Y_p(\F_{p^n}) - p^{2n} - 1$.
For $p=2,3$, let $W_p$ denote the quotient of $H^2_{\et}(\Ybar_p,\Q_l)(1)$
by the image $V_p$ of $\Lambda_p \otimes \Q_l$ under the second homomorphism in
Proposition \ref{boundNS}, and let $\Phi_p$ denote the action of
Frobenius on $W_p$. Since $F_p$ acts trivially on $V_p$, we have
$\Tr \Phi_p^n =\Tr F_p^n-\Tr F_p^n|V_p = \Tr F_p^n -4$ for all $n \geq 0$,
and $f_{F_p}=f_{F_p|V_p}\cdot f_{\Phi_p}=(t-1)^4 f_{\Phi_p}$,
where $f_T$ stands for the
characteristic polynomial of the linear operator $T$. From the traces of
the first $s>0$ powers of a linear operator one can derive the
first $s$ coefficients of its characteristic
polynomial, see \cite{picone}, Lemma 2.4.
Once enough coefficients of $f_{\Phi_p}$ are computed,
the full polynomial $f_{\Phi_p}$ follows from the functional equation
$f_{\Phi_p}(1/x) = \pm x^{-\dim W_p} f_{\Phi_p}(x)$.
Putting all this together, we find $f_{F_p} =
\frac{1}{p}(t-1)^4f_{\Phi_p}$ with
\begin{align*}
f_{\Phi_2} = & 2t^{18} + 2t^{16}+t^{15} + 2t^{14} + t^{13} +
 2t^{12} + t^{11} + 3t^{10} +  \\
   & + 3t^8 + t^7 + 2t^6 + t^5 + 2t^4 + t^3 + 2t^2 + 2, \\
f_{\Phi_3} = & 3t^{18}+5t^{17}+6t^{16}+
         5t^{15}+5t^{14}+6t^{13}-6t^{11}-5t^{10}+ \\
   &  -6t^9-5t^8-6t^7+6t^5+5t^4+5t^3+6t^2+5t+3. \\
\end{align*}
Note that the coefficient of $t^9$ in $f_{\Phi_2}$ is zero, so we
used the number of points over $\F_2^{10}$ to compute the
coefficient of $t^8$, from which we determined
the sign of the functional equation to be positive.
Both $f_{\Phi_p}$ are irreducible. Their roots are not
integral and therefore not roots of unity. By Proposition
\ref{boundNS} we find that the Picard numbers of $\Ybar_2$ and
$\Ybar_3$ are both bounded by $4$, so they are equal to $4$ and
$\Lambda_p$ has finite index in $\NS(\Ybar_p)$ for $p=2,3$. From Lemma
\ref{uptosq} we conclude that up to a square factor the discriminants of
$\NS(\Ybar_2)$ and $\NS(\Ybar_3)$ are equal to $-28$ and $-32$
respectively. From the first injection of Proposition \ref{boundNS} we
find $\rk \NS(\Ybar) \leq 4$. Suppose we had equality. Then the lattice
$\NS(\Ybar)$ would be isomorphic to a sublattice of finite index in
$\NS(\Ybar_p)$ for both $p=2$ and $p=3$. By Lemma \ref{uptosq} this
implies that up to a square factor, the discriminant of $\NS(\Ybar)$
is equal to both $-28$ and $-32$. This contradicts the fact that
$-28$ and $-32$ do not differ by a square factor. We therefore
conclude that equality does not hold and we have $\rk \NS(\Ybar) \leq 3$.
Since the classes $D_1$, $D_2$, and $D_3$ are linearly independent, we
deduce $\rk \NS(\Ybar) =3$.
\endproof

\section{Nonexistence of canonical vector heights}

As in \cite{barorig}, we let $\sigma_i$ denote the involution
associated to the $2$--to--$1$ projection $Y\ra \P^1 \times \P^1$
along the $i$-th copy of $\P^1$ in $\P^1 \times \P^1\times\P^1$, and
for $i,j,k \in \{1,2,3\}$, we set $\sigma_{ijk} = \sigma_i \sigma_j
\sigma_k$. Let $\D^* = \{D_1^*,D_2^*,D_3^*\}$ be the basis that is
dual to the basis $\D = \{D_1,D_2,D_3\}$ of $\NS(\Ybar) \otimes \R$.
Then
$$
\vec h = \sum_{i=1}^3 h_{D_i} D_i^*,
$$
is a vector height, so for every divisor class $E \in \NS(\Ybar)
\otimes \R$, a Weil height $h_E$ associated to $E$ is up to $O(1)$
given by $P \mapsto \vec h(P) \cdot E$.  In our computations, we use
the heights $h_{D_i}$ defined by $\pi_i$
and the usual logarithmic height on $\P^1(\Bbb Q)$.

Suppose $\ss$ is an automorphism of $Y$ and that the pullback
$\ss^*$ acting on $\NS(\Ybar)\otimes \R$ has a real eigenvalue $\ww>1$ with
associated eigenvector $E$. Silverman \cite{Silver} defined the
canonical height (with respect to $\ss$) to be
\[
\hat h_E(P)=\lim_{n\to \infty} \ww^{-n}h_E(\ss^n P).
\]
This height is canonical with respect to the automorphism $\ss$,
since $\hat h_E(\ss P)=\ww \hat h_E(P)$.

Set $\gamma = \frac{1}{2}(1+\sqrt{5})$. Then $\alpha$ and $\omega$
in \cite{barorig} are equal to $\gamma^2$ and $\gamma^6$
respectively. Suppose $(i,j,k)$ is a permutation of $(1,2,3)$. The
eigenvector $E_{ijk}$ of $\sigma_i^*\sigma_j^*\sigma_k^* =
\sigma_{kji}^*$ associated to the eigenvalue $\omega$, as defined
in \cite{barorig}, equals
$\frac{1}{2}\gamma(-D_i+\gamma D_j + \gamma^2 D_k)$. Set $P_0 =
([0:1],[0:1],[0:1])$. Table \ref{hts} contains the estimates
$\omega^{-n}\vec h (\sigma_{kji}^n P_0) \cdot E_{ijk}$ to the
canonical height $\hat{h}_{E_{ijk}}(P_0)$ (canonical with respect to
$\sigma_{kji}$) for all permutations $(i,j,k)$ and $n \in
\{1,\ldots, 5\}$.

\begin{table}[t]
$$
\begin{array}{|c|c|c|c|c|c|c|}
\hline n&(1,2,3) &(1,3,2) &(2,1,3) &(2,3,1) &(3,1,2) &(3,2,1) \cr
\hline
 1 & 0.3438678 & 1.0306631 & 1.7914641 & 2.0624775 & 1.7723601 &1.6340533 \cr
 2 & 0.4711022 & 1.0326396 & 1.8311032 & 2.1288087 & 1.8613679 &1.7950761 \cr
 3 & 0.4745990 & 1.0365615 & 1.8328300 & 2.1330968 & 1.8675712 &1.7982461 \cr
 4 & 0.4747015 & 1.0364020 & 1.8329385 & 2.1332594 & 1.8679417 &1.7986626 \cr
 5 & 0.4746928 & 1.0364196 & 1.8329585 & 2.1332721 & 1.8679467 &1.7986781 \cr
 \hline
\end{array}
$$
\caption{Estimates for $\hat h_{E_{ijk}}(P_0)$ for the permutations
$(i,j,k)$ of $(1,2,3)$.} \label{hts}
\end{table}

These estimates appear to converge geometrically, as expected. We
believe, without rigorous proof, that the estimates of the canonical
heights for $n=5$ are correct up to an error of at most $0.0001$,
and are probably correct up to $0.00001$. The following theorem
therefore gives evidence against the existence of a canonical vector
height on $Y$.

\begin{theorem}
If the estimates $\omega^{-5} h_{E_{ijk}}(\sigma_{kji}^5 P_0)$ in
Table \ref{hts} are equal to the canonical heights
$\hat{h}_{E_{ijk}}(P_0)$ up to an absolute error of at most $\error$, then
the surface $Y$ does not admit a canonical vector height.
\end{theorem}
\begin{proof}
Suppose a canonical vector height $\hat{\vec h}$ exists on $Y$. Then
for every permutation $(i,j,k)$ of $(1,2,3)$ we get a linear
equation (see \cite{barorig})
$$
\hat{\vec h}(P_0) \cdot E_{ijk} = \hat{h}_{E_{ijk}}(P_0).
$$
The three permutations $(3,2,1)$, $(2,3,1)$, and $(3,1,2)$ give
three linearly independent equations from which we can compute the
coefficients $a_i$ in $\hat{\vec h}(P_0) = \sum_{i=1}^3 a_i D_i^*$.
We get
\begin{equation}\label{lineq}
\left[\mymatrix{a_1& a_2& a_3}\right]A =
\left[\mymatrix{\hat{h}_{E_{321}}(P_0)& \hat{h}_{E_{231}}(P_0)&
\hat{h}_{E_{312}}(P_0)}\right],
\end{equation}
where $A$ is the matrix whose columns contain the coefficients with
respect to the basis $\D$ for $E_{321}$, $E_{231}$, and $E_{312}$,
respectively.  That is,
\[
A=\frac{\cc}{2}\left[\begin{matrix}
                  \cc^2 & \cc^2 & \cc \\
                  \cc   & -1    & \cc^2 \\
                  -1    & \cc   & -1
  \end{matrix}\right].
\]
The absolute values of the entries of $A^{-1}$ are bounded by
$2\gamma^{-1}$, so when we use the estimates of
$\hat{h}_{E_{ijk}}(P_0)$ for $n=5$ in Table \ref{hts}, the solution
$$
(a_1, a_2 , a_3) = (0.719498, 0.805119, 0.963093)
$$
to (\ref{lineq}) is accurate up to $\varepsilon =
3(2\gamma^{-1})(\error)$. From $E_{123} =
\frac{1}{2}\gamma(-D_1+\gamma D_2 +\gamma^2 D_3)$ we find that, up
to an absolute error of at most $\frac{1}{2}\gamma(1+\gamma
+\gamma^2)\varepsilon \approx 1.571$, the canonical height
$\hat{h}_{E_{123}}(P_0)$ equals
$$
(0.719498\,D_1^*+ 0.805119\,D_2^*+ 0.963093\,D_3^*)\cdot E_{123}
\approx 2.51169.
$$
This contradicts the estimate in the first column of Table
\ref{hts}, so we conclude that $Y$ does not admit a canonical vector
height.
\end{proof}

\end{document}